\documentclass{article}
\usepackage{amsmath,amssymb,amsbsy,amsthm}
\let\emptyset\varnothing
\let\epsilon\varepsilon
\let\phi\varphi

\def\N{\mathbb N}

\def\S{\mathcal S}

\def\cl{\operatorname{cl}}

\newtheorem{theorem}{Theorem}
\newtheorem{lemma}{Lemma}

\newtheorem{definition}{Definition}

\begin{document}

\title{A criterion for hypothesis testing for stationary processes}
\author{Daniil Ryabko \\ {\em INRIA Lille-Nord Europe,}\\ {\em 40, Avenue Halley
59650 Villeneuve d'Ascq, France} \\ {\tt daniil@ryabko.net}}
\date{}
\maketitle
\begin{abstract}
Given a finite-valued sample $X_1,\dots,X_n$ we wish to test whether it was generated by a stationary ergodic process 
belonging to a family $H_0$, or it was generated by a stationary ergodic process outside $H_0$. 
We require the Type I error of the test to be  uniformly bounded, while the type II error 
has to be mande not more than a finite number of times with probability 1. 
For this notion of consistency we provide necessary and sufficient conditions on the family $H_0$
for the existence of a consistent test. 
This criterion is illustrated with applications to testing for a  membership to parametric families,
generalizing some existing results. In addition, we analyze a stronger notion of consistency, which requires
finite-sample guarantees on error of both types, and provide some necessary and some sufficient conditions for 
the existence of a consistent test.
We emphasize that no assumption on the process distributions are made beyond stationarity and ergodicity.
\end{abstract}
\centerline{\em Keywords: Hypothesis testing, stationary processes, ergodic processes, distributional distance.}

\section{Introduction}
Given a sample $X_1,\dots,X_n$ (where $X_i$ are from a finite  alphabet $A$) that is known to be generated by a stationary ergodic process,
 we wish to decide whether it was generated by a distribution  belonging to a  certain family $H_0$, versus it was generated by a stationary 
ergodic  distribution  that does not belong to $H_0$. 
Unlike most of the works on the subject, we do not 
assume that $X_i$ are i.i.d., but only make a much weaker assumption that
the distribution generating the sample is stationary ergodic. 

A test is a function that takes a sample and an additional parameter $\alpha$ (the significance level), and  gives a  binary (possibly incorrect) answer:  the sample
was generated by a distribution from $H_0$ or by a stationary ergodic distribution not belonging to $H_0$. 
Here we are concerned with characterizing those families   $H_0$  for which consistent tests exist.

We consider the following notion of consistency.
Call a test {\em  consistent} if,  for any pre-specified {\em level} $\alpha\in(0,1)$, any sample size $n$ and 
any distribution in $H_0$ {\em the probability of Type I error (the test says ``not $H_0$'')
is not greater than $\alpha$}, while for every stationary ergodic distribution  from outside $H_0$ and every $\alpha$ {\em  Type II error 
(the test says $H_0$) is made only a finite number of times (as the sample size goes to infinity) with probability 1}. 
This notion of consistency represents  a classical statistical
approach to the problem, and suites well situations where the hypothesis $H_0$ is considerably more simple than the alternative,
for example when $H_0$ consists of just one distribution, or when it is some parametric family, or when 
it is the hypothesis of homogeneity or that of independence. 

\noindent {\bf Prior work.} 
There is a vast body of literature on hypothesis testing for i.i.d. (real- or discrete-valued) data (see e.g. \cite{leh}).
In the context of discrete-valued  i.i.d. data, the necessary and sufficient conditions 
for the existence of a consistent test  are rather  simple to obtain: there is a consistent test for $H_0$ (against ``i.i.d. but not $H_0$'') if and only if  $H_0$
is closed,
where the topology is that of the parameter space (probabilities of each symbol), e.g. see  \cite{cs}.  
The consistency being easy to ensure, the prime concern for the case of i.i.d. data  is optimality. 

There is, however,  much less literature on hypothesis testing beyond i.i.d. or parametric models, while the questions
of determining whether a consistent test exists (for different notions of consistency and different hypotheses) is much less trivial.
For a weaker notion of consistency, namely, requiring that the test should stabilize on the correct answer for a.e. realization of the process (under either $H_0$ or 
the alternative), 
 \cite{kief} constructs a consistent test for so-called
constrained finite-state model classes (including finite-state Markov and hidden Markov processes), against the general alternative of stationary ergodic processes. 
For the same  notion of consistency, 
 \cite{nob} gives sufficient conditions on two families $H_0$ and $H_1$ that consist of stationary ergodic real-valued processes, under which a consistent continuous
test exists, extending the results of~\cite{dp} for i.i.d. data. The latter condition is that $H_0$ and $H_1$ are contained in disjoint $F_\sigma$ sets (countable unions of closed sets), with respect to the topology of weak convergence.
For the notion of consistency that we consider, consistent tests for some specific hypotheses, but under the general alternative of stationary ergodic
processes, have been proposed in 
  \cite{ra,rg,rr}, which address problems of testing identity, independence, estimating the order of a Markov process,
and also the change point problem. Some impossibility results for testing hypotheses about stationary ergodic processes can be found in \cite{mor2, dr09}.

\noindent  {\bf The results.}
The aim of this work is to provide topological characterizations of the hypotheses for which consistent
tests exist, for the case of stationary ergodic distributions. The obtained characterization is rather
similar to those mentioned above for the case of i.i.d. data, but is with respect to the topology 
of distributional distance (or weak convergence). The fact that necessary and sufficient conditions are obtained  
 indicates that this topology is  the right one to consider.

A distributional distance between two process distributions  is defined as a
weighted sum of probabilities of all possible tuples $X\in A^*$, where $A$ is the alphabet and  the weights
are positive and have a finite sum. 
The main result is the following theorem (formalized in the next sections).

\noindent{\bf Theorem.} There exists a consistent test for $H_0$ 
 if and only if
$H_0$ has probability 1 with respect to ergodic decomposition of every distribution from the closure  of $H_0$.  

The test  that we construct to establish this result  is based on empirical estimates of distributional distance.
For a given level $\alpha$,  it  takes the largest $\epsilon$-neighbourhood
of the closure of $H_0$ that has probability not greater than $1-\alpha$ with respect to every ergodic process in it, 
and  outputs $0$ if the sample falls into this neighbourhood, and $1$ otherwise.

To illustrate the applicability of the main result,  we show that the families of  $k$-order Markov processes and
 $k$-state Hidden Markov processes (for any natural $k$), satisfy the conditions of the theorem, and therefore there exists a consistent test 
for membership to these families.

It should be emphasized that the results of this work concern what is possible in principle; finding an
efficient testing procedure for each specific hypothesis for which we can demonstrate 
existence of a consistent test is a different problem.

\section{Preliminaries}\label{s:prel}
Let $A$ be a finite alphabet, and denote  $A^*$ the set of words (or tuples) $\cup_{i=1}^\infty A^i$ and $A^\infty$ the set of all one-way infinite sequences.
For a word  $B\in A^*$ the symbol $|B|$ 
stands for the length of $B$.
 Distributions, or  (stochastic) processes, are measures on the space $(A^\infty,\mathcal F_{A^\infty})$, where $\mathcal F_{A^\infty}$ is the 
Borel sigma-algebra of $A^\infty$.
%
Denote $\#(X,B)$ the number of occurrences of a word $B\in A^*$ in a
word $X\in A^*$ and $\nu(X,B)$ its frequency:
$$
\#(X,B)=\sum_{i=1}^{|X|-|B|+1} I_{\{(X_i,\dots,X_{i+|B|-1})=B\}},
$$
and
\begin{equation}\label{eq:freq}
\nu(X,B)=\left\{\begin{array}{cc}{1\over|X|-|B|+1}\#(X,B) & \text{ if }|X|\ge|B|,\\ 0 &\text{
otherwise,}\end{array}\right.
\end{equation}
where $X=(X_1,\dots,X_{|X|})$.
For example, $\nu(0001,00)= 2/3.$

We use the abbreviation $X_{1..k}$ for $X_1,\dots,X_k$.
A process $\rho$ is {\em stationary}  if 
$$
\rho(X_{1..|B|}=B)=\rho(X_{t..t+|B|-1}=B)
$$
 for any $B\in A^*$ and $t\in\N$.
Denote  $\S$ the set of all stationary  processes on $A^\infty$.
A stationary process $\rho$ is called {\em (stationary) ergodic} if
   the frequency of occurrence of each word
$B$ in a sequence $X_1,X_2,\dots$ generated by $\rho$ tends to its
a priori (or limiting) probability a.s.: 
$\rho(\lim_{n\rightarrow\infty}\nu(X_{1..n},B)= \rho(X_{1..|B|}=B))=1$. 
By virtue of the ergodic theorem (e.g. \cite{bil}), this definition can be shown to be equivalent 
to the standard definition of stationary ergodic processes (every shift-invariant set has measure 0 or 1; see e.g. \cite{cs}). 
 Denote $\mathcal E$ the set  of all stationary ergodic processes.

\begin{definition}[distributional distance]
 The  distributional distance  is defined for a pair of processes
$\rho_1,\rho_2$ as follows~\cite{gray}:
$$
d(\rho_1,\rho_2)=\sum_{k=1}^\infty w_k |\rho_1(X_{1..|B_k|}=B_k)-\rho_2(X_{1..|B_k|}=B_k)|,
$$
where $w_k=2^{-k}$ and  $B_k$, $k\in\N$  range through the set $A^*$ of all words in length-lexicographical order 
(the   weights and ordering are fixed for the sake of concreteness only).
\end{definition}
It is easy to see that $d$ is a metric.
Equipped with this metric, the space of all stochastic processes is separable  and complete; moreover, it's a compact. 
The set  of stationary 
processes  $\S$ is its convex closed subset  (hence a  compact too). 
The set of all finite-memory stationary distributions is dense in $\mathcal S$. (Taking only those that have rational transition probabilities
we obtain a countable dense subset of $\S$.)
The set $\mathcal E$ is not convex (a mixture of stationary ergodic distributions is always stationary but never ergodic) and is  not closed (its closure is $\mathcal S$).
We refer to \cite{gray} for more details and proofs of these facts.

When talking about closed and open subsets of $\S$ we assume the topology of~$d$.
Compactness of the set $\S$ is one of the main ingredients in the proofs of the main results.
Another is that the distance $d$ can be consistently estimated, as the next lemma shows.

Considering 
the Borel (with respect to the metric $d$) sigma-algebra $\mathcal F_\S$ on the set $\S$, we obtain
a standard probability space $(\S,\mathcal F_\S)$.
An important tool that will be used in the analysis is  {\bf ergodic decomposition} of stationary processes  (see e.g. \cite{gray, bil}):
 which we recall here.  Any stationary process can be expressed as a mixture of stationary ergodic
processes;
more formally, for any $\rho\in\S$ there is a measure $W_\rho$ on $(\S,\mathcal F_\S)$, such 
that $W_\rho(\mathcal E)=1$, and $\rho(B)=\int d W_\rho(\mu)\mu(B)$, for any $B\in\mathcal F_{A^\infty}$.
The {\em support} of a stationary distribution $\rho$ is the minimal closed set $U\subset\S$ such that $W_\rho(U)=1$.

A {\bf test} is a function $\psi^\alpha: A^*\rightarrow\{0,1\}$ that takes as input a sample and a parameter $\alpha\in(0,1)$, and outputs a binary answer,
where the answer $0$ is interpreted as ``the sample was generated by a distribution that belongs to $H_0$'', and the answer 1 as 
``the sample was generated by a stationary ergodic distribution that does not belong to $H_0$.''
A test $\phi$ makes the {\em Type I} error if it says $1$ while $H_0$ is true,
and it makes {\em Type II} error if it	 says $0$ while $H_0$ is false.

\begin{definition}[consistency]
 Call a  test $\psi^\alpha, \alpha\in(0,1)$    consistent as a test of $H_0$ against $H_1$ if:
\begin{itemize}
\item[(i)] The probability of Type I error is always bounded by $\alpha$: $\rho\{X\in A^n: \psi^\alpha(X)=1\}\le\alpha$ for every $\rho\in H_0$,
every $n\in\N$ and every $\alpha\in(0,1)$, and
\item[(ii)]  Type II error is made not more than a finite number of times with probability 1: 
 $\rho(\lim_{n\rightarrow\infty} \psi^\alpha(X_{1..n})=1)=1$ 
for every  
$\rho\in H_1$ and every $\alpha\in(0,1)$. 
\end{itemize}
\end{definition}

\section{Main results}\label{s:main}
The test constructed below is based on {\em empirical estimates of  the distributional distance $d$}:
$$
\hat d(X_{1..n},\rho)=\sum_{i=1}^\infty w_i |\nu(X_{1..n},B_i)-\rho(B_i)|,
$$
where $n\in\N$, $\rho\in\S$, $X_{1..n}\in A^n$. That is, $\hat d(X_{1..n},\rho)$ measures
the discrepancy between empirically estimated and theoretical probabilities.
 For a sample $X_{1..n}\in A^n$ and a hypothesis $H\subset\mathcal E$ define 
$$
\hat d(X_{1..n},H)=\inf_{\rho\in H} \hat d(X_{1..n},\rho).
$$

Construct the  test  $\psi_{H_0}^\alpha, \alpha\in(0,1)$ as follows.
For each  $n\in\N$, $\delta>0$ and $H\subset\mathcal E$ define the neighbourhood  $b^n_\delta(H)$  of $n$-tuples around $H$ as 
  $$b^n_\delta(H):=\{X\in A^n: \hat d(X,H)\le\delta\}.$$
Moreover, let 
$$
\gamma_n(H,\theta):= \inf \{\delta: \inf_{\rho\in H} \rho(b^n_\delta(H))\ge\theta\}
$$ be the smallest radius of a neighbourhood around $H$ that has probability not less than $\theta$ with respect to every process in $H$, and let $C^n(H,\theta):=b^n_{\gamma_n(H,\theta)}(H)$ be a neighbourhood of this radius.
Define  
$$
\psi^\alpha_{H_0}(X_{1..n}):=\left\{\begin{array}{ll} 0& \text{ if }  X_{1..n}\in C^n(\cl H_0\cap\mathcal E,1-\alpha), \\ 1& \text{ otherwise.} \end{array}\right.
$$
We will often omit the subscript $H_0$ from $\psi^\alpha_{H_0}$ when it can cause no confusion.

The main result of this work is the following theorem, whose proof is given in section~\ref{s:proofs}.
\begin{theorem}\label{th:asym}
 Let $H_0\subset\mathcal E$. The following  statements are equivalent:
\begin{itemize}
 \item[(i)] There exists a consistent test for $H_0$ against $\mathcal E\backslash H_0$.
 \item[(ii)] The test $\psi^\alpha_{H_0}$ is   consistent. 
 \item[(iii)] The set $H_0$ has probability 1 with respect to ergodic decomposition of every $\rho$ in the closure of $H_0$: 
              $W_\rho(H_0)=1$ for each $\rho\in\cl H_0$.
\end{itemize}
\end{theorem}

\section{Examples}\label{s:ex}
The first simple illustration of Theorem~\ref{th:asym} above is identity testing, or goodness of fit: testing whether a distribution 
generating the sample obeys a certain given law, versus it does not. Let $\rho\in\mathcal E$, $H_0=\{\rho\}$. 
Since $H_0$ is closed, Theorem~\ref{th:asym} implies that there is a consistent test for $H_0$. 
Identity testing is a classical problem of mathematical statistics, with solutions (e.g. based on Pearson's $\chi^2$ statistic) for i.i.d.
data (e.g. \cite{leh}), and Markov chains \cite{bilmark}. 
 For stationary ergodic processes, \cite{rg} gives a consistent test when $H_0$ has a finite and bounded memory, and \cite{rr} for the general case.

Another example is bounding the order of  a Markov or  a Hidden Markov process.
Theorem~\ref{th:asym} implies that for any given $k\in\N$ there is a consistent  test 
of the hypothesis $\mathcal M^k$= ``the process is Markov of order not greater than $k$'' (against $\mathcal E\backslash \mathcal M^k$).
Moreover,  there is a  consistent test of $\mathcal{HM}^k$=``the process is given by a Hidden Markov process with not more than $k$ states.''
Indeed, in both cases ($k$-order Markov, Hidden Markov with not more than $k$ states), the hypothesis $H_0$
is a parametric family, with a compact set of parameters, and a continuous function mapping parameters 
to processes (that is, to the space $\mathcal S$). Weierstrass theorem then implies that the image of such a compact
parameter set is closed (and compact). Moreover, in both cases $H_0$ is closed under taking ergodic decompositions.
Thus, by Theorem~\ref{th:asym}, there exists a consistent test.

The problem of estimating the order of a (hidden) Markov process, based on a sample from it, was addressed in a number of works. In the contest 
of hypothesis testing,  consistent tests for $\mathcal M^k$ against $\mathcal M^t$ with $t>k$ were given in \cite{andgoo}, see also \cite{bilmark}.
For a weaker notion of consistency (the test has to stabilize on the correct answer eventually, with probability 1) 
the existence of a consistent test for $\mathcal{HM}^k$   was established in \cite{kief}. 
For the notion of consistency considered here, a consistent 
test for $\mathcal M^k$  was proposed in \cite{ra}, while 
for the case of testing $\mathcal {HM}^k$  the result above is apparently new.

\section{Uniform testing}
Finally, let us consider a stronger notion of hypothesis testing, that requires uniform speed of convergence for
errors of either type.

A test $\phi$ is called {\bf uniformly consistent} if for every $\alpha$ there is an $n_\alpha\in\N$ such that
for every $n\ge n_\alpha$ the probability of error on a sample of size $n$ is less than $\alpha$: $\rho(X\in A^n: \phi(X)=i) <\alpha$ for every $\rho\in H_{1-i}$ and every $i\in\{0,1\}$.

For $H_0,H_1\subset\S$, the {\bf uniform test} $\phi_{H_0,H_1}$ is constructed   as follows. For each $n\in\N$ let
\begin{equation}
 \phi_{H_0,H_1}(X_{1..n})\\:=\left\{\begin{array}{ll} 0& \text{ if }  \hat d(X_{1..n},\cl H_0\cap\mathcal E)<\hat d(X_{1..n},\cl H_1\cap\mathcal E), \\ 1 & \text{ otherwise.} \end{array}\right.
\end{equation}

\begin{theorem}[uniform testing]\label{th:uni}
 Let $H_0\subset\S$ and $H_1\subset\S$.  If 
$W_\rho(H_i)=1$ for every $\rho\in \cl H_i$
then 
the test $\phi_{H_0,H_1}$ is uniformly consistent.
Conversely, if there exists a uniformly consistent test for $H_0$ against $H_1$ then 
$W_\rho(H_{1-i})=0$ for any $\rho\in cl H_i$.
\end{theorem}
The proof is given in the next section.

\section{Proofs}\label{s:proofs}
The proof of the main results will use the following lemmas.
\begin{lemma}[$\hat d$ is consistent]\label{th:dd} Let $\rho,\xi\in\mathcal E$ and let a sample $X_{1..k})$
be generated by $\rho$.  Then
$$
\lim_{k\rightarrow\infty}\hat d(X_{1..k},\xi)=d(\rho,\xi)\ \rho\text{-a.s.}
$$
\end{lemma}
The proof is based on the fact that the frequency of each word converges to its expectation.
For each $\delta$  we can find a time by which the first  $K(\delta)$  frequencies will have converged up to $\delta$,
where $K(\delta)$ is such that  the cumulative weight of the rest of the frequencies is smaller than $\delta$ too. 

\begin{proof}
For any $\epsilon>0$ find such an index $J$ that
$\sum_{i=J}^\infty w_i<\epsilon/2$.
For each $j$ we have $\lim_{k\rightarrow\infty}\nu(X_{1..k},B_j)= \rho(B_j)$
a.s., so that
$|\nu(X_{1..k},B_j) - \rho(B_j)|<\epsilon/(2Jw_j)$ from some $k$ on; denote  $K_j$ this $k$.
 Let
$K=\max_{j<J}K_j$ ($K$ depends
on the realization $X_1,X_2,\dots$).
 Thus, for $k>K$ we have
\begin{multline*}
   |\hat d(X_{1..k},\xi) - d(\rho,\xi)| 
= \left|\sum_{i=1}^\infty w_i\big(|\nu(X_{1..k},B_i)-\xi(B_i)| - |\rho(B_i)-\xi(B_i)| \big) \right|
  \\
   \le \sum_{i=1}^\infty w_i|\nu(X_{1..k},B_i)-\rho(B_i)|  
    \le \sum_{i=1}^J w_i|\nu(X_{1..k},B_i)-\rho_X(B_i)|  +\epsilon/2 
   \\
   \le \sum_{i=1}^Jw_i \epsilon/(2Jw_i) +\epsilon/2 = \epsilon,
\end{multline*}
which proves the statement.
\end{proof}

\begin{lemma}[smooth probabilities of deviation]\label{th:count}
 Let $m>2k>1$, $\rho\in\S$, $H\subset\S$, and $\epsilon>0$. Then
\begin{equation}\label{eq:count1} 
\rho(\hat d(X_{1..m},H)\ge\epsilon)  
  \le 2\epsilon'^{-1} \rho(\hat d(X_{1..k},H)\ge \epsilon'),
\end{equation}
where $\epsilon':=\epsilon - \frac{2k}{m-k+1} - t_k$ with $t_k$ being the sum of  all the weights of tuples longer than $k$ in the definition of $d$: $t_k:=\sum_{i:|B_i|>k}w_i$. Further,
\begin{equation}\label{eq:count2} 
\rho(\hat d(X_{1..m}, H)\le\epsilon)
  \le 2\rho\left(\hat d(X_{1..k},H)\le \frac{m}{m-k+1}2\epsilon + \frac{4k}{m-k+1}\right).
\end{equation}
\end{lemma}
The meaning of this lemma is as follows.  For any word $X_{1..m}$, if it is far away from (or close to)
a given distribution $\mu$ (in the empirical distributional distance), then some of its shorter subwords $X_{i..i+k}$ are far from (close to) $\mu$ too.
 In other words, for a stationary distribution $\mu$,  it cannot happen that  a small 
 sample is likely to be close to $\mu$, but a larger sample is likely to be far.

\begin{proof}
Let $B$ be  a tuple such that $|B|<k$ and $X_{1..m}\in A^m$ be any sample of size $m>1$. 
The number of occurrences of $B$ in $X$ can be bounded by the number of occurrences of $B$ in subwords of $X$ of length $k$ as follows:
\begin{multline*}
\#(X_{1..m},B) \le   \frac{1}{k-|B|+1}\sum_{i=1}^{m-k+1}\#(X_{i..i+k-1},B) +  2k \\ =\sum_{i=1}^{m-k+1}\nu(X_{i..i+k-1},B)+  2k.
\end{multline*}
Indeed, summing over $i=1..m-k$ the number of occurrences of $B$ in all $X_{i..i+k-1}$ we count each
occurrence of $B$ exactly $k-|B|+1$ times, except for those that occur in the first and last $k$ symbols.
Dividing by $m-|B|+1$, and using the definition (\ref{eq:freq}), 
 we obtain
\begin{equation}\label{eq:ff}
\nu(X_{1..m},B)  \le \frac{1}{m-|B|+1}\left(\sum_{i=1}^{m-k+1}\nu(X_{i..i+k-1},B)| + 2k\right).
\end{equation}
Summing over all $B$, for any $\mu$, we get
\begin{equation}\label{eq:fracs}
\hat d(X_{1..m},\mu)  \le \frac{1}{m-k+1} \sum_{i=1}^{m-k+1}\hat d(X_{i..i+n-1},\mu) + \frac{2k}{m-k+1}+t_k,
\end{equation}
where in the right-hand side $t_k$ corresponds to all the summands in the left-hand side for which $|B|>k$, where for the rest of the summands we used $|B|\le k$.
Since this holds for any $\mu$, 
we conclude that
\begin{multline}\label{eq:ne}
\hat d(X_{1..m}, H)  \le \frac{1}{m-k+1} \left(\sum_{i=1}^{m-k+1}\hat d(X_{i..i+k-1}, H)\right)  + \frac{2k}{m-k+1}+t_k.
\end{multline}
Note that the $\hat d(X_{i..i+k-1}, H)\in[0,1]$. Therefore, for the average in the r.h.s.\ of~\eqref{eq:ne} to be 
larger than $\epsilon'$, at least $\epsilon'/2 (m-k+1)$ summands have to be larger than $\epsilon'/2$.

Using stationarity, we can conclude
$$
\rho\left(\hat d(X_{1..k},H)\ge\epsilon'\right)\ge \epsilon'/2 \rho\left(\hat d(X_{1..m}, H)\ge\epsilon\right),
$$
proving~(\ref{eq:count1}). The second statement can be proven similarly; indeed, analogously to~(\ref{eq:ff}) we have
 \begin{multline*}
\nu(X_{1..m},B)  \ge \frac{1}{m-|B|+1}\sum_{i=1}^{m-k+1}\nu(X_{i..i+k-1},B)-  \frac{2k}{m-|B|+1} \\ \ge \frac{1}{m-k+1}\left(\frac{m-k+1}{m}\sum_{i=1}^{m-k+1}\nu(X_{i..i+k-1},B)\right) -  \frac{2k}{m},
\end{multline*} where we have used $|B|\ge 1$. Summing over different $B$,  we obtain (similar to~(\ref{eq:fracs})),
 \begin{equation}\label{eq:ne2}
\hat d(X_{1..m},\mu)  \ge \frac{1}{m-k+1} \sum_{i=1}^{m-k+1}\frac{m-k+1}{m}\hat d_k(X_{i..i+n-1},\mu) - \frac{2k}{m}
\end{equation}(since the frequencies are non-negative, there is no $t_n$ term here).
For the average in~\eqref{eq:ne2} to be smaller than $\epsilon$, at least half of the summands must be smaller than $2\epsilon$.
Using stationarity of $\rho$, this implies~(\ref{eq:count2}).
\end{proof}

\begin{lemma}\label{th:cont} Let $\rho_k\in\S$, $k\in\N$ be a sequence of processes that converges to a process $\rho_*$.
Then, for any $T\in A^*$ and $\epsilon>0$ if $\rho_k(T)>\epsilon$ for infinitely many indices $k$, then
$\rho_*(T)\ge\epsilon$
\end{lemma}
\begin{proof}
 The statement follows from the fact that  $\rho(T)$ is continuous as a function of $\rho$. 
\end{proof}

\smallskip

\noindent{\em Proof of Theorem~\ref{th:asym}.} 
The implication {\em (ii) $\Rightarrow$ (i)} is obvious. We will show {\em (iii) $\Rightarrow$ (ii)} and {\em (i) $\Rightarrow$ (iii)}.
To establish the former, we have to show that the family of tests $\psi^\alpha$ is consistent.
By construction, for any $\rho\in\cl H_0\cap\mathcal E$ we have $\rho(\psi^\alpha(X_{1..n})=1)\le\alpha$. 

To prove the consistency of $\psi$, it remains to show that $\xi(\psi^\alpha(X_{1..n})=0)\rightarrow0$ a.s. for any $\xi\in\mathcal E\backslash H_0$ and $\alpha>0$.
To do this, fix any $\xi\in \mathcal E\backslash H_0$ and let $\Delta:=d(\xi,\cl H_0):=\inf_{\rho\in \cl H_0\cap\mathcal E} d(\xi,\rho)$. Since $\cl H_0$
is closed, we have $\Delta>0$. Suppose that there exists an $\alpha>0$, such that, for infinitely many $n$, some samples from the $\Delta/2$-neighbourhood of $n$-samples
around $\xi$ are sorted as $H_0$ by $\psi$, that is, $C^n(\cl H_0\cap\mathcal E, 1-\alpha)\cap b_{\Delta/2}^n(\xi)\ne\emptyset$. Then
for these $n$ we have $\gamma_n(\cl H_0\cap\mathcal E,1-\alpha)\ge\Delta/2$. 

This means that there exists an increasing sequence $n_m,m\in\N$, and a 
sequence $\rho_m\in \cl H_0$, $m\in\N$, such that 
  $$\rho_{m}(\hat d(X_{1..n_{m}},\cl H_0\cap\mathcal E)> \Delta/2)>\alpha.$$
Using Lemma~\ref{th:count}, (\ref{eq:count1}) (with $\rho=\rho_{m}$, $m=n_m$, $k=n_k$, and $H=\cl H_0$),  and taking  $k$  large enough to have $t_{n_k}<\Delta/4$,
 for every $m$  large enough to have $\frac{2n_k}{n_m-n_k+1}<\Delta/4$, we obtain 
\begin{equation}\label{eq:down-}
8\Delta^{-1}\rho_{m}\left(\hat d(X_{1..n_k},\cl H_0)  \ge\Delta/4\right)  \ge  \rho_{m}\left(\hat d(X_{1..n_m},\cl H_0)\ge\Delta/2\right) > \alpha. 
\end{equation}
%
%
Thus, 
\begin{equation}\label{eq:new}
 \rho_{m}(b^{n_k}_{\Delta/4}(\cl H_0\cap\mathcal E))<1-\alpha\Delta/8.
\end{equation}

Since the set $\cl H_0$ is compact (as a closed subset of a compact set $\S$),  we may assume (passing to 
a subsequence, if necessary) that  $\rho_m$  converges
to a certain $\rho_*\in\cl H_0$.
Since~\eqref{eq:new} this holds for infinitely many $m$, 
using Lemma~\ref{th:cont} (with $T=b^{n_k}_{\Delta/4}(\cl H_0\cap\mathcal E)$) we conclude that
 $$\rho_*(b^{n_k}_{\Delta/4}(\cl H_0\cap\mathcal E))\le1-\Delta\alpha/8.$$
Since the latter inequality holds for infinitely many indices $k$ we also have 
$$
\rho_*(\limsup_{n\rightarrow\infty}\hat d(X_{1..n},\cl H_0\cap\mathcal E)>\Delta/4)>0.
$$
However,  we must have $\rho_*(\lim_{n\rightarrow\infty}\hat d(X_{1..n},\cl H_0\cap\mathcal E)=0)=1$
for every $\rho_*\in\cl H_0$: indeed, for $\rho_*\in\cl H_0\cap\mathcal E$ it follows from Lemma~\ref{th:dd}, and for $\rho_*\in\cl H_0\backslash\mathcal E$ 
from  Lemma~\ref{th:dd}, ergodic decomposition and the conditions of the theorem ($W_\rho(H_0)=1$ for $\rho\in\cl H_0$).

This contradiction shows that for every $\alpha$ there are not more than finitely many $n$ for which $C^n(\cl H_0\cap\mathcal E, 1-\alpha)\cap b_{\Delta/2}^n(\xi)\ne\emptyset$. 
To finish the proof of the  implication,  it remains to note that, as follows from Lemma~\ref{th:dd}, 
\begin{multline*}
\xi\{X_1,X_2,\dots.:X_{1..n}\in b_{\Delta/2}^n(\xi)\text{ from some $n$ on}\}\\\ge \xi\left(\lim_{n\rightarrow\infty}\hat d(X_{1..n},\xi)=0\right) =1.
\end{multline*}

To establish the implication {\em (i) $\Rightarrow$ (iii)},
 we assume that there exists a  
 consistent test $\phi$ for $H_0$, and we will show that
 $W_\rho(\mathcal E\backslash H_0)=0$ for every $\rho\in \cl H_0$. 
Take $\rho\in \cl H_0$ and suppose that  $W_\rho(\mathcal E\backslash H_0)=\delta>0$. 
We have 
\begin{multline*}
 \limsup_{n\to\infty} \int_{\mathcal E\backslash H_0} d W_\rho(\mu) \mu(\psi^{\delta/2}_n=0)  \le  \int_{\mathcal E\backslash H_0} \limsup_{n\to\infty} d W_\rho(\mu) \mu(\psi^{\delta/2}_n=0)=0,
\end{multline*}
where the inequality follows from Fatou's lemma (the functions under integral are all bounded by 1), and the equality from the consistency of $\psi$. 
Thus, from some $n$ on 
we will have $\int_{\mathcal E\backslash H_0} d W_\rho \mu(\psi^{\delta/2}_n=0) < 1/4$ so that $\rho(\psi^{\delta/2}_n=0)<1-3\delta/4$. 
For any set $T\in A^n$ the function $\mu(T)$ is continuous as a function of $\mu$. In particular, it holds for the
set $T:=\{X_{1..n}:\psi_n^{\delta/2}(X_{1..n})=0\}$. Therefore,
since $\rho\in\cl H_0$, for any $n$ large enough we can find a $\rho'\in H_0$ such that $\rho'(\psi^{\delta/2}_n=0)<1-3\delta/4$,
 which contradicts
the consistency of $\psi$. Thus, $W_\rho(H_0)=1$, and Theorem~\ref{th:asym} is proven.
\qed

\noindent{\em Proof of Theorem~\ref{th:uni}.}
To prove the first statement of the theorem,
we will show that the test $\phi_{H_0,H_1}$ is a uniformly consistent test for $\cl H_0\cap\mathcal E$ against $\cl H_1\cap\mathcal E$ (and hence for $H_0$ against $H_1$),
under  the conditions of the theorem. 
Suppose that, on the contrary,  for some $\alpha>0$ for every $n'\in\N$ there is a process $\rho\in \cl H_0$ such that 
$\rho(\phi(X_{1..n})=1)>\alpha$ for some $n>n'$. 
Define 
$$
\Delta:=d(\cl H_0,\cl H_1):=\inf_{\rho_0\in \cl H_0\cap\mathcal E, \rho_1\in \cl H_1\cap\mathcal E} d(\rho_0,\rho_1),
$$ which is positive since $\cl H_0$ and $\cl H_1$ are closed and disjoint.
We have
\begin{multline}\label{eq:union}
\alpha<\rho(\phi(X_{1..n})=1)\\ \le   \rho(\hat d(X_{1..n},H_0)\ge\Delta/2\ or \ \hat d(X_{1..n},H_1)<\Delta/2)\\ \le \rho(\hat d(X_{1..n},H_0)\ge\Delta/2) + \rho(\hat d(X_{1..n},H_1)<\Delta/2).
\end{multline} 
This implies that either $\rho(\hat d(X_{1..n},\cl H_0)\ge\Delta/2)>\alpha/2$ or $\rho(\hat d(X_{1..n},\cl H_1)<\Delta/2)>\alpha/2$,
so that, by assumption, at least one of these inequalities holds for infinitely many $n\in\N$ for some sequence  $\rho_n\in H_0$.
Suppose that it is the first one, that is, there is an increasing sequence $n_i$, $i\in\N$ and a sequence $\rho_i\in\cl H_0$, $i\in\N$ 
such that 
\begin{equation}\label{eq:da}
\rho_i(\hat d(X_{1..n_i},\cl H_0)\ge\Delta/2)>\alpha/2 \text{ for all }i\in\N. 
\end{equation}
The set $\S$ is compact, hence so is its closed subset $\cl H_0$. Therefore, the sequence $\rho_i$, $i\in\N$ must
contain a subsequence that converges to a certain process $\rho_*\in\cl H_0$. Passing to a subsequence if necessary, we may assume
that this convergent subsequence is the sequence $\rho_i$, $i\in\N$ itself.

Using Lemma~\ref{th:count}, (\ref{eq:count1}) (with $\rho=\rho_{n_m}$, $m=n_m$, $k=n_k$, and $H=\cl H_0$),  and taking  $k$  large enough to have $t_{n_k}<\Delta/4$,
 for every $m$  large enough to have $\frac{2n_k}{n_m-n_k+1}<\Delta/4$, we obtain 
\begin{equation}\label{eq:down}
8\Delta^{-1}\rho_{n_m}\left(\hat d(X_{1..n_k},\cl H_0)  \ge\Delta/4\right)  \ge  \rho_{n_m}\left(\hat d(X_{1..n_m},\cl H_0)\ge\Delta/2\right) > \alpha/2. 
\end{equation}
That is, we have shown that for any large enough index $n_k$ the inequality 
$\rho_{n_m}(\hat d(X_{1..n_k},\cl H_0)\ge\Delta/4)> \Delta\alpha/16$ holds for 
infinitely many indices $n_m$. 
From this and Lemma~\ref{th:cont} with  $T=T_k:=\{X:\hat d(X_{1..n_k},\cl H_0)\ge\Delta/4\}$  we conclude 
that $\rho_*(T_k)>\Delta\alpha/16$.
The latter holds for infinitely many $k$; that is,  $\rho_*(\hat d(X_{1..n_k},\cl H_0)\ge\Delta/4)>\Delta\alpha/16$ infinitely often.
Therefore, 
$$
\rho_*(\limsup_{n\rightarrow\infty} d(X_{1..n},\cl H_0)\ge\Delta/4)>0.
$$ However, we must have 
$$
\rho_*(\lim_{n\rightarrow\infty} d(X_{1..n},\cl H_0)=0)=1
$$
for every $\rho_*\in\cl H_0$: indeed, for $\rho_*\in\cl H_0\cap\mathcal E$ it follows from Lemma~\ref{th:dd}, and for $\rho_*\in\cl H_0\backslash\mathcal E$ 
from  Lemma~\ref{th:dd}, ergodic decomposition and the conditions of the theorem. 

Thus, we have arrived at a contradiction that shows that $\rho_n(\hat d(X_{1..n},\cl H_0)>\Delta/2)>\alpha/2$ cannot hold for 
infinitely many $n\in\N$ for any sequence of  $\rho_n\in\cl H_0$. Analogously, we can show that $\rho_n(\hat d(X_{1..n},\cl H_1)<\Delta/2)>\alpha/2$
cannot hold for infinitely many $n\in\N$ for any sequence of $\rho_n\in\cl H_0$. Indeed,  using Lemma~\ref{th:count}, equation~(\ref{eq:count2}), we can show that
  $\rho_{n_m}(\hat d(X_{1..n_m},\cl H_1)\le\Delta/2) > \alpha/2$ for a large enough $n_m$
implies $\rho_{n_m}(\hat d(X_{1..n_k},\cl H_1)\le 3\Delta/4)> \alpha/4$ for a smaller $n_k$.
Therefore, if we assume that $\rho_n(\hat d(X_{1..n},\cl H_1)<\Delta/2)>\alpha/4$ for infinitely many $n\in\N$ for some sequence of $\rho_n\in\cl H_0$, then 
we will also find a $\rho_*$ for which $\rho_*(\hat d(X_{1..n},\cl H_1)\le 3\Delta/4)> \alpha/4$ for infinitely
many $n$, which, using Lemma~\ref{th:dd} and ergodic decomposition, can be shown to contradict the fact that $\rho_*(\lim_{n\rightarrow\infty}  d(X_{1..n},\cl H_1)\ge\Delta)=1$.
 
Thus, returning to~(\ref{eq:union}), we have shown that from some $n$ on there is no $\rho\in \cl H_0$ for which $\rho(\phi=1)>\alpha$
holds true. The statement for $\rho\in\cl H_1$ can be proven analogously, thereby finishing the proof of the first statement.

To prove the second statement of the theorem, 
we assume that there exists a uniformly consistent test $\phi$ for $H_0$ against $H_1$, and we will show that
$W_\rho(H_{1-i})=0$ for every $\rho\in \cl H_i$.
Indeed, let $\rho\in \cl H_0$, that is, suppose that there is a sequence $\xi_i\in H_0, i\in\N$ such that $\xi_i\to\rho$. 
Assume  $W_\rho(H_1)=\delta>0$ and take $\alpha:=\delta/2$.
Since the test $\phi$ is uniformly consistent, there is an $N\in\N$ such  that  for every $n>N$ we have  
\begin{multline*}
 \rho(\phi(X_{1..n}=0))\le \int_{H_1} \phi(X_{1..n}=0)dW_\rho + \int_{\mathcal E\backslash H_1} \phi(X_{1..n}=0) dW_\rho \\ \le \delta\alpha + 1-\delta \le 1-\delta/2.
\end{multline*}
Recall that, for $T\in A^*$,  $\mu(T)$ is a continuous function in $\mu$. In particular, this holds for the set $T=\{X\in A^n: \phi(X)=0\}$, 
for any given $n\in\N$. Therefore, for every  $n>N$ and for every  $i$ large enough,  $\rho_i(\phi(X_{1..n})=0)<1-\delta/2$ implies also $\xi_i(\phi(X_{1..n})=0)<1-\delta/2$
which contradicts $\xi_i\in H_0$.
This contradiction shows $W_\rho(H_1)=0$ for every $\rho\in \cl H_0$. The case $\rho\in \cl H_1$ is analogous.
\qed

\end{document}